\documentclass[10pt,a4paper]{article}
\usepackage[utf8]{inputenc} 
\usepackage{amsmath} 
\usepackage{amsfonts}
\usepackage{amssymb}
 
 \usepackage{cleveref}
 
\usepackage{epsfig}
\usepackage{graphics}
 
\newtheorem{Theorem}{Theorem}

\newtheorem{Remark}[Theorem]{Remark}

\newtheorem{Assumption}[Theorem]{Assumption}

\newcommand{\zaa}{\alpha}
 
\newcommand{\zg}{\gamma}
\newcommand{\ZDE}{\delta}

\newcommand{\zl}{\lambda} 
 
 \newcommand{\ZOMq}{\Omega} 
 \newcommand{\zthe}{\theta}
\newcommand{\zt}{\tau}

\newcommand{\zben}{\beta_n}

\newcommand{\zthen}{\theta_n}
\def\zn{z_n}
\def\xin{\xi_n}
\def\zln{\zl_n}

\def\zsuNO{\sum _{n=1}^{+\ZIN}}

\newcommand{\zzr}{\mathbb{R}}

\newcommand{\intt}{\int_0^t}
\newcommand{\ints}{\int_0^s}
\newcommand{\intr}{\int_0^r}

 \newcommand{\zdia}{~~\rule{1mm}{2mm}\par\medskip}
 
\newcommand{\ZIN}{\infty}

 \newcommand{\ZD}{\;\mbox{\rm d}}
 \newcommand{\ZSUno}{\sum_{n=1}^{\ZIN}}

 \newcommand{\ZLA}{\label} 
\newcommand{\ZBI}{\bibitem}
  
  \def\intmint{\int_{-\ZIN}^t}
  \def\intzpin{\int_0^{+\ZIN}}
  \def\intL{\int_0^L}
 
\author{
L. Pandolfi\thanks{Dipartimento di Scienze Matematiche ``Giuseppe Luigi Lagrange'', Politecnico di Torino, Corso Duca degli Abruzzi 24, 10129 Torino, Italy (luciano.pandolfi@polito.it)}
}
\title{Joint identification via deconvolution of the flux and energy relaxation kernels of the Gurtin-Pipkin model in thermodynamics with memory\thanks{
This papers fits into the research program of the GNAMPA-INDAM and has been written in the framework of the   ``Groupement de Recherche en Contr\^ole des EDP entre la France et l'Italie (CONEDP-CNRS)''.}}

\begin{document}

\maketitle

\begin{abstract}
In this paper we present a linear method for the identification of both the energy and flux relaxation kernels in the equation of thermodynamics with memory proposed by M.E. Gurtin and A.G. Pipkin. The method reduces the identification of the two kernels to the solution of two (linear) deconvolution problems. The energy relaxation kernel is  reconstructed by means of energy measurements as the solution of a Volterra integral equation of the first kind which does not depend on the still unknown flux relaxation kernel.   Then, flux measurements are used to identify the flux relaxation kernel.
\end{abstract}

\noindent{\bf Key Words} Relaxation kernel, diffusion processes with memory, identification, deconvolution
\section{Introduction}

The linearized version of the model proposed in~\cite{GurtinPIPKIN} by Gurtin and Pipkin to describe thermodynamical processes with memory depends on two memory kernels $\beta(t)$, the \emph{energy relaxation kernel,} and $a(t)$, the \emph{flux  relaxation kernel.} These kernels are material  properties of the body and have to be identified using suitable measurements. Several identification methods have been proposed   in engineering and mathematical papers  (see the references in \cref{Sect-References}). In most of the cases, these methods assume $\beta=0$. Here we are going  to extend to the case $\beta\neq 0$ a linear algorithm first proposed in~\cite{PandINVprob} for the identification of the   kernel $a(t)$ when $ \beta(t)=0 $. The bonus of this algorithm is that it reduce the identification of the kernels to   (linear) deconvolution problems.

In order to understand the identification algorithm, we shall shortly present the derivation of the Gurtin-Pipkin model in  \cref{sect:THEmodel}. We shall see that the evolution in time of the temperature is described by the following equation
\begin{equation}
\ZLA{eqGP}
\frac{\ZD}{\ZD t}\left [\zthe(t)+\intmint \beta(t-s)\zthe(s)\ZD s\right ]=\intmint a(t-s)\Delta\zthe(s)\ZD s \,.
\end{equation}
 
Here $\zthe=\zthe(x,t)$ is a function of the time $t$ and of a variable $x$ in a region $\ZOMq$   ($\Delta $ is the laplacian in the variable $x$).

Note that the variable $x$ and also the time variable will not be explicitly indicated unless needed for clarity.

As explained in~\cite{BykovMILLISECONDS}, the identification of the relaxation kernels can be obtained using a sample in the form of a bar.
 Hence, we shall consider the case $\ZOMq=(0,L)$ so that $\Delta\zthe =\zthe_{xx}$. But, the extension of the algorithm to general regions is important for applications to nondestructive testing. This extension  is reserved to the future.

We state the assumptions which we use to justify the identification algorithm:
\begin{Assumption}
\begin{enumerate}
\item $a(t)$ and $\beta(t)$ are bounded and  of class $C^2(0,+\ZIN)$ (the derivatives have a continuous extension to $t=0$).
\item $a(0)>0$  (This condition is particularly important since it implies that signals propagate  with finite velocity  and for this reason \cref{eqGP} is also called the \emph{hyperbolic heat equation}).

\end{enumerate}   
\end{Assumption}
In fact, the derivation of \cref{eqGP} assumes also that $ a(t) $ and $ \beta(t)$ and their derivatives    are integrable on $ [0,+\ZIN) $ (and so $a(t)$, $\beta(t)$    tend to zero for $ t\to+\ZIN $).
Thermodynamics   impose stronger conditions to the kernels (for example, it is proved in~\cite{DayDECREASkernels} that
 under natural conditions the kernels decreases    and certain ``positivity'' conditions are proved in~\cite{AmendolaLIBRO}). These properties are not used in  the identification algorithm, but they can be used to improve the numerical implementation, as discussed in~\cite{PandINVprob}. The properties of the kernels $a(t)$ and $\beta(t)$  (which are not used in the justification of our algorithm) imply that the system dissipates energy.
 
The organization of the paper is as follows: in \cref{sect:THEmodel} we describe informally the identification algorithm. In order to understand the rational under the algorithm, we describe also the derivation of \cref{eqGP}.

\Cref{sec:preliSUeq}  presents preliminary results on the solutions of \cref{eqGP} (which of course depend on suitable initial and boundary conditions) while the algorithm is justified in \cref{sec:justiALGO}  while   \cref{Sect-References} contains comments to the literature.

\section{\ZLA{sect:THEmodel} Informal description of the identification algorithm}

In order to understand the rationale behind the identification algorithm, we must shortly describe the derivation of \cref{eqGP} in~\cite{GurtinPIPKIN}. Note that \cref{eqGP} is a linearized version of a ``true'' nonlinear model (see also~\cite{brandon} for a nice description of the nonlinear model) and so $ \zthe $ does not represent the absolute temperature but it is the perturbation of a stationary temperature of a nonlinear process. Hence, $\zthe  $ can be zero.

The derivation of \cref{eqGP} depends on the following two principles:
\begin{enumerate}
\item the (density  of the) flux of heat $q(x,t)$ at position $x$ and time $t$ depends on the past history of the gradient of the temperature $\zthe(x,t)$ at the same position $x$, and it is given by
\begin{equation}
\ZLA{eq:delFLUSSO}
q(x,t)=-\intzpin a(s)\nabla\zthe(x,t-s)\ZD s=-\intmint a(t-s)\nabla\zthe(x, s)\ZD s
\end{equation}
 \Cref{eq:delFLUSSO} replace the usual Fourier Law $q(x,t)=-\nabla\zthe(x,t)$.
 \item the internal energy is related to the temperature by the relation
 \begin{multline}
e(x,t)
=b+c\zthe(x,t)+\intzpin \beta(s)\zthe(x,t-s)\ZD s\\
= b+c\zthe (x,t) +\intmint \beta(t-s)\zthe(s)\ZD s \,.
 \end{multline}
  This relation replace the usual relation $e(x,t)=b+c\zthe (x,t)$. In the following we assume that the scale has been chosen so to have $c=1$.
\end{enumerate}

\Cref{eqGP}   follows when we balance the energy at every time $t$  and position $x$, i.e. we impose
\[
\frac{\ZD}{\ZD t} e(x,t)=-\nabla\cdot q(x,t)\,.
\]

 We note that \cref{eqGP}  can be written in several equivalent forms, in particular
 \begin{equation}
 \ZLA{eqGP0}
 \zthe'(t)+\beta(0)\zthe(t)+\intmint \beta'(t-s)\zthe(s)\ZD s=\intmint a(t-s)\zthe_{xx}(s)\ZD s
 \end{equation}
(we recall: $c=1$).

We did not yet specify the initial and boundary conditions. We fix an initial time $t_0$, it is not restrictive to put $t_0=0$, and we assume that the measurements are taken for $t>0$. So, in order to solve
\cref{eqGP0} we need the initial condition
\begin{equation}\ZLA{iniCONDI}
\zthe(x,s)=\xi_0(x,s)\quad \mbox{given for $s<0$}\,,\quad \zthe(x,0)=\xi(x)\,.
\end{equation} 
The function $ \xi_0(x,s)$ needs not be continuous so that in general $\xi_0(x,0)$ does not even makes sense (and it is not equal to $\xi(x)$). Discontinuity at zero can be practically realized by putting the body in contact with a suitable source. This observation is used in the identification algorithm since  we assume $ \zthe(x,s)=\xi_0(x,s)=0 $ $ s<0 $ and $\xi$ possibly different from zero.
 
 When $\zthe(x,t)=\xi_0(x,t)=0$ for $t<0$,  \cref{eqGP0} takes the form
 
\begin{multline}
 \ZLA{eqGP1}
 \frac{\ZD}{\ZD t}\left [\zthe(t)+\intt \beta(t-s)\zthe(s)\ZD s\right ]\\
 = 
  \zthe'(t)+\beta(0)\zthe(t)+\intt \beta'(t-s)\zthe(s)\ZD s=\intt a(t-s)\zthe_{xx}(s)\ZD s \,.
 \end{multline}

Due to the fact that $a(t)$ is related to the flux, it is clear that its identification will require the measure of the flux (which can be measured only at the ends of the bar).  The kernel $\beta(t)$ is related to the relaxation of energy, and the quantity that can be directly measured is the temperature. So, in order to identify $\beta(t)$ we measure the temperature.
In the case of a bar, conceivably we can measure the temperature at each point.

The relaxation kernels are identified by taking the following measurements of the sample of the material:
\begin{enumerate}
\item first we impose boundary conditions $ \zthe(0,t)=\zthe(L,t)=0 $ and an initial temperature  $ \zthe(x,0)=\xi(x) $. We measure
\begin{enumerate}
\item\ZLA{eitemBETA1} the total enery as a function of time. But, the quantity which can be measured is the temperature. So we estimate
\[ 
\Theta_\xi(t)=\int_0^L \zthe(x,t)\ZD x\,;
 \]
 \item
 \ZLA{eitemA1}
  the flux at $ x=L $, i.e.
 \[ Y_\xi(t)
 =q(L,t)=-\int_0^t a(t-s)\zthe_x(L,s)\ZD s=-\lim _{ x\to L^- }
 \int_0^t a(t-s)\zthe_x(x,s)\ZD s\,.
  \]
\end{enumerate}
\item then we repeat the same measurements but with $ \xi=0 $ while $ \zthe(0,t)=f(t) $ and $ \zthe(L,t)=0 $:
\begin{enumerate}
\item\ZLA{eitemBETA2}the total enery as a function of time. In fact we measure
\[ 
\Theta^f(t)=\int_0^L \zthe(x,t)\ZD x\,;
 \]
 \item
 \ZLA{eitemA2}
  the flux at $ x=L $, i.e.
 \[ Y^f(t)
 =q(L,t)=-\int_0^t a(t-s)\zthe_x(L,s)\ZD s=-\lim _{ x\to L^- }
 \int_0^t a(t-s)\zthe_x(x,s)\ZD s\,.
  \]
\end{enumerate}
\end{enumerate}

\begin{Remark}{\rm
Note the use of a nonzero initial temperature. While there is no problem to apply a time varying temperature $f(t)$ at one end of the bar,  the initial condition is more difficult to realize. We shall see  that we need a special initial condition which is easily realized.\zdia
}\end{Remark}

The energy relaxation kernel $ \beta(t) $ can be computed from the measurements described in the \cref{eitemBETA1,eitemBETA2} even if $ a(t) $ is still unknown. 

Once $ \beta(t) $ has been computed then $ a(t) $ is computed from the measurements described in the~\cref{eitemA1,eitemA2}.

\begin{Remark}\ZLA{RemaSULtempeFRONTIE}{\rm
We observe that the integral $\int_0^L\zthe(x,t)\ZD x $ can be estimated since in practice the kernels can be determined from samples which have the form of  a bar. It has an interest to see that identification of $\beta(t)$ can be acieved also by measuring the temperature $\zthe(L,t)$ when 
we impose the boundary condition $ \zthe_x(L,t)=0 $ (i.e. the right hand side is now insulated) and we measure
\begin{description}
\item[i)]\ZLA{eitemBETA1BIS} (instead of \cref{eitemBETA1}) we measure $ \zthe(L,t) $ when
\[ 
\zthe(x,0)=\xi(x)\,, \quad \zthe(0,t)=0\,,\quad  \zthe_x(L,t)=0 \,.
 \]
 
\item[ii)]\ZLA{eitemBETA2BIS} (instead of \cref{eitemBETA2})
we measure $ \zthe(L,t) $ when
\[ 
\zthe(x,0)=0\,,\qquad \zthe(0,t)=f(t)\,,\qquad  \zthe_x(L,t)=0 \,. 
 \]  
 
\end{description}
 This variant is examined in \cref{sect:variants}.\zdia
 }
 \end{Remark}

\subsection{\ZLA{Sect-References}References to previous results}
  \Cref{eqGP} has important applications in thermodynamics and viscoelasticity and in other applications (like nonfickian diffusion, which is often encountered in biological matherials). Hence the determination of the parameters in the equation has been widely studied in 
  mathematics and in engineering journal. See the reference in~\cite{PandINVprob}. The methods used in Engineering journals are mostly based on this idea: the kernels to be identified are assumed to belong to specific classes, for example Prony sums, and depend on few parameters. The equation is   discretized and the parameters of the kernels are determined so that the discretized version of the solution best fit the experimental measure $y(t)$. Usually $y(t)$ is a measure of temperature or flux or, in viscoelasticity, it is a measure of the traction on the boundary. This method was first proposed to identify one memory kernel (see for example~\cite{BykovMILLISECONDS}) but of course it can be used as well to identify the
   two memory kernels in \cref{eqGP}. In this contest we mention in particular the paper~\cite{GolubMaslovFERNATI} which is concerned wit three dimensioanl viscoelasticity. In this case, different relaxation kernels correspond to the shear and bulk  creep relaxation kernels.
   
   The method proposed in~\cite{LorenziSINESTRARI} is as follows. It is noted that  the observation $y(t)$ 
is a linear    functional of $\zthe$, which is unknown since \cref{eqGP} cannot be solved, because the kernels are still unknown. But then, the equation of $y$ and  \cref{eqGP} constitute a nonlinear system  in the unknown  $(\zthe, N)$. The solution of this system provides both $\zthe$ and the unknown kernel $N$. Usually this system is solved by first computing the derivative of $y$ so to get a nonlinear differential equations in the unknowns $(\zthe, N)$  in a suitable Hilbert spaces. The method, first proposed in the case that the equation depends on only one kernel was then extended to the identification of two
    kernels, for example in~\cite{ColomboGuidetti,JannoLOrenzi,LoreROCCA,PAISjanno}. We note that this method requires the solution of highly nonlinear equation in a Hilbert space but it is noted in~\cite{Guidetti} (in the case of one unknown kernel, but the observation is easily extendable) that once the kernels and $\zthe$ have been identified in a first interval $[0,\zt]$ then, thanks to the Volterra structure of the problem, the identification up to $T$ is reduced to a linear problem in a Hilbert space. In contrast with this, the method in this paper requires the solution of two Volterra integral equations of the first kind in $\zzr$, since the quantity to be identified are the scalar functions $\beta(t)$ and $a(t)$.

\section{\ZLA{sec:preliSUeq}Representation of the solutions}
The solutions of \cref{eqGP1}  have been studied by several authors. A nice reference is~\cite{Belleni}. Here we need a Fourier type approach which is similar to the one used in~\cite{PANDlibro}. We consider the operator $A$ in $L^2(0,L)$ which is defined as follows:
  
\[ 
{\rm dom}\,A=H^1_0(0,L)\cap H^2(0,L)\,,\qquad A\zthe=\zthe_{xx}\,.
 \]
 Let \[ 
 \phi_n(x)=\zg_0\sin  \zl_n x\,,\qquad \zg_0=\sqrt{\frac{2}{L}
}\,,\quad \zln=n\frac{\pi	}{L}\,.  \]

The function $ \phi_n $ is an eigenvector of $ A $
whose eigenvalue is $ -\zln^2 $  
and
$ \{\phi_n\}   $ is an orthonormal basis  $ L^2(0,L) $.
 We expand the solutions of \cref{eqGP1} with conditions
 
\[
\zthe(0)=\xi
  \quad {\rm and}  \quad \left\{\begin{array}{l}
  \zthe(0,t)=f(t)\,,\\
  \zthe(L,t)=0  \,,\\
\end{array}\right.
\]
Let
\[
\left\{\begin{array}{l}
\xin=\intL \xi(x) \phi_n(x)\ZD x\,,\\
\zthen(t)=\intL \zthe(x,t)\phi_n(x)\ZD x\,, 
\end{array}\right.
\qquad \mbox{so that}\qquad 
\left\{\begin{array}{l}
\xi(x)=\ZSUno\xin\phi_n(x)\,,\\
\zthe(t)=\ZSUno \phi_n(x)\zthen(t)\,.
\end{array}\right.
\]
The function $\zthen(t)$ solves
 \begin{multline}
 \ZLA{eqDIztheN}
\frac{\ZD}{\ZD t}\left [
 \zthen(t)+\intt \beta(t-s)\zthen(s)\ZD s
 \right ]= \zthen'(t)+\beta(0)\zthen(t)+\intt \beta'(t-s)\zthen(s)\ZD s
 \\=-\zl_n^2\intt a(t-s)\zthen(s)+\phi_n'(0)\intt a(t-s)f(s)\ZD s\,,\qquad \zthen(0)=\xin\,.
 \end{multline}
Note that
\[ 
 \frac{\phi_n'(0)	 }{\zl_n}=\zg_0=\sqrt{\frac{2}{L}} \,.
 \]

We intoduce the functions $\zn(t)$ which solve
 \begin{multline}
 \ZLA{eqDIzETaN}
\frac{\ZD}{\ZD t}\left [
 \zn(t)+\intt \beta(t-s)\zn(s)\ZD s
 \right ]=\zn'(t)+\beta(0)\zn(t)+\intt \beta '(t-s)\zn(s)\ZD s 
 \\=-\zln^2\intt a(t-s)\zn(s)\,,\qquad  \zn(0)=1\,.
 \end{multline}
We note that
\[
\zthen(t)=\zn(t)\xin-\intt\frac{\ZD}{\ZD s}\left [\zthen(t-s)\zn(s)\right ]\ZD s\,.
\]
 We compute explicitly the derivative under the integral sign and we see that
 \begin{equation}
 \ZLA{eq:soluDItHETAnante}
 \zthen(t)= \zn(t)\xin+\phi_n'(0)\intt \zn(t-s) \left [
 \ints a(s-r)f(r)\ZD r
 \right ]\ZD s\,.
 \end{equation}

 The last integral can be manipulated as follows (when $f(t)$ is differentiable):
 \begin{multline}\ZLA{eq:ManiULTint}
\phi_n'(0) \intt \zn(t-s) \left [
 \ints a(s-r)f(r)\ZD r
 \right ]\ZD s\\ 
 =-\frac{\phi_n'(0)}{\zln^2}\intt f(t-s)\left [-\zl_n^2\ints a(s-r)\zn(r)\ZD r\right ]\ZD s\\
= -\frac{\zg_0}{\zln}\intt f(t-s)\left [
 \frac{\ZD}{\ZD s}\left (
 \zn(s)+\ints \beta(s-r)\zn(r)\ZD r
 \right)
 \right ]\ZD s \\
= -\frac{\zg_0}{\zl_n}\left [f (0)\left (
 \zn(t)+\intt \beta(t-r)\zn(r)\ZD r
 \right)
 -f(t)\right.\\
 \left.+\intt f'(t-s)\left (
 \zn(s)+\ints \beta(s-r)\zn(r)\ZD r
 \right)\ZD s
 \right ]\,.
 \end{multline}
In the following, we shall use $f\in C^1$ with $f(0)=0$ (consistent with $ \zthe(x,0)=0 $), i.e.
\begin{equation}\ZLA{eq:expreF}
f(t)=\intt g(s)\ZD s\,.
\end{equation}
Hence  
\begin{multline}
\ZLA{eq:soluDItHETAn}
\zthe(t)=\ZSUno \phi_n(x)\zthe_n(t)\,,\qquad 
 \zthen(t)=\zn(t) \xin +\frac{\zg_0}{\zl_n}\intt g(s)\ZD s\\
 -\frac{\zg_0}{\zl_n} \intt g(t-s)\left (
 \zn(s)+\ints \beta(s-r)\zn(r)\ZD r
 \right)\ZD s\,.
\end{multline}
 
\section{\ZLA{sec:justiALGO}Identification of the kernels}

Here we show that the proposed method enables us to identify the relaxation kernels by solving (linear) deconvolution problems. 
 
First we consider the identification of the energy relaxation kernel $ \beta(t) $.

\paragraph{Identification of the energy relaxation kernel.}
Note that 
\[ 
\int_0^L\phi_n(x) \ZD x=
 \zaa_n\frac{1}{\zl_n}\,,\qquad \zaa_n=
 \left\{\begin{array}{cl}
0&\mbox{if $ n $ is even}\\
2\zg_0&\mbox{if $ n $ is odd.}
\end{array}\right.
 \]  
Then it is easily computed that,  when $ f(t)=0 $:

\begin{equation}\ZLA{serie1}
\int_0^L\zthe(x,t)\ZD x= \ZSUno \frac{1}{\zl_n}\zaa_n\xi_n z_n(t)\,.
 \end{equation}
 We use the special initial condition $\xi(x)$ whose Fourier coefficients are 
 $ \xi_n=1/\zl_n $  The function provided  by this measurement is denote $ H(t) $:
\[ 
H(t)= \ZSUno \frac{1}{\zl_n^2}\zaa_n  z_n(t)\,.
 \]
 
 \begin{Remark}{\rm
 We shall see that samples with this special initial condition are easily realized in practice without using the eigenvalues and eigenvectors of the problem.\zdia
 }\end{Remark}
 
Now we use  formula~(\ref{eq:soluDItHETAn}) with $ \xi_n=0 $ to compute $ \int_0^L\zthe(x,t)\ZD x $ when $ \zthe(0,t)=f(t) $. We find that
\begin{multline*}
\int_0^L\zthe(x,t)\ZD x= \left (\intt g(s)\ZD s\right ) \ZSUno\zaa_n\frac{1}{\zl_n^2}- \intt g(t-s)H(s)\ZD s\\
- \intt \beta(s)\left [\int_0^{t-s}g(t-s-r)H(r)\ZD r\right ]\,\ZD s\,.
\end{multline*}
The left hand side is measured, $ H(t) $ has been already estimated and so the first and second terms on the right side are known. Hence, $ \beta(t) $ is computed from the deconvolution problem:
\[ 
  \intt \beta(s)\left [\int_0^{t-s}g(t-s-r)H(r)\ZD r\right ]\,\ZD s=\mbox{known function}\,.
 \]

Note that $ a(t) $, still unidentified, does not appear in the previous computations.
  
  \paragraph{\ZLA{sec:IndendDIa}Identification of the flux relaxation kernel.}
  
  It is not possible to measure the flux inside a body. So, in this case we are forced to measure the flux at $L$, when the boundary conditions are 
  \[ 
\zthe(x,0)=\xi(x)\,,\quad   \zthe(0,t)=f(t)\,,\qquad \zthe(L,t)=0\,.
   \]
   and we measure the two fluxes produced either by $ \xi $ or by $ f $.
   
   Let us first consider the flux when $ f=0 $ and $ \xi (x)\neq 0 $. In this case we have
   \begin{align*}
   \zthe(x,t)=\ZSUno \phi_n(x)z_n(t)\xi_n\,,\quad \intt a(t-s)\zthe(x,s)\ZD s=\ZSUno \phi_n(x)\xi_n \intt a(t-s)z_n(s)\ZD s\\
   =-\ZSUno \phi_n(x)\frac{\xi_n}{\zl_n^2}\left (-
   \zl_n^2\intt a(t-s)z_n(s)\ZD s
   \right )\\
   -\ZSUno \phi_n(x)\frac{\xi_n}{\zl_n^2}\frac{\ZD}{\ZD t}\left (z_n(t)+\intt \beta(t-s) z_n(s)\ZD s\right)\,.
   \end{align*}
   Using \[\frac{\ZD}{\ZD x} \phi_n(L)=\zg_0(-1)^n\zl_n \]
   we see that the flux is
   \[ 
   \intt a(t-s)\zthe_x(L,s)\ZD s=\frac{\ZD}{\ZD t}\left [
   \ZSUno (-1)^n\frac{\zg_0\xi_n }{\zl_n}\left (  z_n(t)+\intt \zben(t-s) z_n(s)\ZD s\right )\right ]\,.
    \]
    This flux can be measured, provided that the initial condition $ \xi  $ can be realized.  As
    already stated and  discussed below, we can realize that initial condition $ \xi $ whose Fourier coefficients are
    $ \xi_n=1/\zl_n $. So, we measure the corresponding flux, that we denote $ K(t) $:
    
    \[ 
    K(t)=
    \frac{\ZD}{\ZD t}\left [
   \ZSUno (-1)^n\frac{\zg_0  }{\zl^2_n}\left (  z_n(t)+\intt \beta(t-s) z_n(s)\ZD s\right )\right ]\,.
     \]
     Now we compute the flux throughot the right hand $ L$ due to the boundary temperature $ \zthe(0,t)=f(t) $ ($ f(t) $ as 
     in \cref{eq:expreF}) and $ \xi =0 $.  
      Using \cref{eq:soluDItHETAn} we see that 
      \begin{multline*}
      \intt a(t-s)\zthe(x,s)\ZD s= \ZSUno \phi_n(x)\frac{\zg_0}{\zl_n}\intt a(s)\int_0^{t-s} g(r)\ZD r\,\ZD s\\
      - \ZSUno \phi_n(x)\frac{\zg_0}{\zl_n }  \intt \left (\intr a(r-s)z_n(s)\ZD s\right )h(t-r)\ZD r
      \end{multline*}
      where
      \[ 
      h(t)=g(t)+\intt \beta(t-r)g(r)\ZD r
       \]
       is a known function  since $ \beta(t) $ has already been identified.
       So we have
       \begin{multline}\ZLA{fluxFpasso2}
       \intt a(t-s)\zthe(x,s)\ZD s=\ZSUno \phi_n(x) \frac{\zg_0}{\zl_n}\intt a(s)\int_0^{t-s}g(r)\ZD r\,\ZD s\\
       +\ZSUno \frac{\zg_0}{\zl_n^3}\phi_n(x)\intt h(t-r)\left (-\zl_n^2\intr a(r-s) z_n(s)\ZD s\right )\ZD r\\
=\ZSUno\phi_n(x)\frac{\zg_0}{\zl_n}\intt a(s)\int_0^{t-s}g(r)\ZD r\,\ZD s\\
+\ZSUno \frac{\zg_0}{\zln^3} \phi_n(x)\intt h(t-r)\frac{\ZD}{\ZD r}\left [
z_n(r)+\intr \beta(r-s)z_n(s)\ZD s
\right ]\ZD r       
       \,.
       \end{multline}
       
       In order to compute the flux throughout $ L $ we must compute the derivative respect to $ x $. We note that $ \frac{\ZD}{\ZD x}\phi_n (L)=(-1)^n\zg_0\zl_n$. So, the contribution of the second line is
\begin{multline}\ZLA{eq:scambioDeriSERIE} 
\intt h(t-r)\frac{\ZD}{\ZD r}\left (
\ZSUno (-1)^n\frac{\zg_0^2}{\zl_n^2}\left [z_n(r)+\intr \beta(r-s)z_n(s)\right ]\ZD s 
\right )\ZD r\\
=\zg_0\intt h(t-r)K(r)\ZD r
 \end{multline}
 and this is a known quantity.
 
 Now we consider the first addendum on the right side of \cref{fluxFpasso2}. We note that
 \[ 
\frac{\ZD}{\ZD x}\ZSUno \phi_n(x) \frac{\zg_0}{\zl_n}
=
 \zg_0\ZSUno \cos n\frac{\pi}{L}x
 =\zg_0\left (\frac{l}{\pi}\delta-\frac{1}{2}\right ) 
  \]
 ($\ZDE$ is the Dirac's delta).  
      So, the flux at $ L $ is
      \[
      Y^f(t)=\frac{\zg_0}{2} \intt  a(s)\int_0^{t-s}g(r)\ZD r\,\ZD s -
   \intt h(t-s) K(s)\ZD s\,.
       \]
     The relaxation kernel $ a(t) $ can be computed from this equality via deconvolution.
     
We observe that the exchange of the series and the derivative in \cref{eq:scambioDeriSERIE}  is easily justified since the sequence of continuous functions $\{z_n(t)\} $ is bounded on bounded intervals.     
     
     This ends the description of the identification procedure. It remains to understand whether it is possible to realize easily  a source with the special temperature 
     \[
\xi_0=\ZSUno \frac{1}{\zl_n}\phi_n(x)=\zg_0\ZSUno\frac{1}{\zl_n} \sin n\frac{\pi}{L} x    \,. 
     \]
 The positive answer was given in~\cite{PandAPPLmath}. The temperature $ \xi (x) $ is the solution of the following problem:
 
      \[
      \xi_{xx}=0\,,\qquad \xi(0)=1\,,\quad \xi(L)=0\,.
       \]
  So, it is the stationary temperature achieved by a thermal body whose boundary temperatures are 
  \[
   \zthe(0,t)=0\,,\qquad \zthe(L,t)=\frac{1}{\zg_0}\,. 
  \]
  In fact, it is known that the heat equation $\eta'=\eta_{xx}$ with conditions $\eta(0,t)=\eta(L,t)=0$ is exponentially stable, so that $\lim _{t\to+\ZIN}\eta(x,t)=0$ in $L^2(0,L)$ for every initial condition $\eta(x,0)=\xi(x)\in L^2(0,L)$.

      Let   $\zthe$ solves 
      \begin{equation}\ZLA{eq:StandHEAT}
      \zthe'(t)=\Delta\zthe\,,\qquad \zthe(0)=0 \quad \zthe(0,t)=0\,,\quad \zthe(L,t)=\zg_0\,. 
\end{equation}
The function
\[
\eta(x,t)=\zthe(x,t)-\xi (x) 
\]
solves the heat equation with the boundary condition put equal zero and the initial condition equal to $\xi $.

Stability implies $\lim_{t \to+\ZIN} \eta(\cdot,t)=0 $ in $L^2(\ZOMq)$ so that
\[
\lim_{t \to+\ZIN}\zthe(x,t)=\xi_0(x)\,:
\]
in order to realize a source with  the special temperature $\xi_0$ needed in the identification process it is sufficient to apply the constant  boundary temperatures at $x=0$, $x=L$ as described in~(\ref{eq:StandHEAT}), for a time large enough.

Note that the same arguments apply not only to the solutions of the standard heat equation but also to the solutions of \cref{eqGP1}, provided that the kernels are so chosen that the system  is stable, as usually happens in practice.
 
     \begin{Remark}{\rm
     The Fourier expansion has been used to justify the procedure, but the actual implementation of the procedure does not need any Fourier expansion, not even in order to realize the special initial condition needed in this procedure.\zdia
     }\end{Remark}
     \subsection{\ZLA{sect:variants}The variant in Remark~\ref{RemaSULtempeFRONTIE}}
We show that the procedure in Remark~\ref{RemaSULtempeFRONTIE} leads to the identification of $ \beta $. We note that in this case the operator $ A $ has to be replaced by
     \[ 
     A=\Delta\,,\qquad {\rm dom}\,A=\{\phi\in H^2(0,L)\,,\  \phi(0)=0\,,\ \phi_x(L)=0\}\,.
      \]
      The eigenvectors of this operator   are
      \begin{equation}\ZLA{eq:neweigenv} 
      \phi_n(x)=\zg_0\sin\left (n+\frac{1}{2}\right )\frac{\pi}{L}x\,,\qquad \zg_0=\sqrt{\frac{2}{L}}
       \end{equation}
       with eigenvalues $ -\zl_n^2 $ where now
      $
        \zl_n=\left (n+\frac{1}{2}\right )(\pi/L)
        $.
        
        Note that 
        \[
\phi_n(L)=\zg_0(-1)^n\,.        
        \]
        Granted this  new meanings of the symbols, the formulas for $ \zthe_n(x,t) $, $ z_n(t) $   are still given by~\cref{eqDIzETaN,eq:soluDItHETAnante}. The manipulations in~(\ref{eq:ManiULTint}) still hold and $ \zthe(x,t) $ has the expansion $ \zthe(t)=\ZSUno \zthe_n(t)\phi_n(x) $ in 
        \cref{eq:soluDItHETAn}  where now  $ \phi_n $ is in~(\ref{eq:neweigenv}). So, the algorithm can be justified as follows.

We consider the case $ f=0 $ and an initial temperature $ \xi (x)\neq 0  $ imposed at time $ t=0 $. In this case 
\[ 
\zthe(t,L)=\zg_0\zsuNO  (-1)^n\xi_n z_n(t) \,.
 \]
 Convergence of this series is easily seen because $\{\xi_n\}\in l^2$ and the same estimates as in~\cite{PandIEOT} shows that
 \[
z_n(t)=\cos\zl_n t+\frac{1}{\zl_n}M_n(t) 
 \]
 where $\{M_n(t)\}$ is a sequence of continuous functions which is bounded on bounded intervals.
 
We call $ H(t) $ the temperature $ \zthe(L,t) $ when
\[ 
\xi_0(x)=\zsuNO \frac{1	}{\zl_n} \phi_n(x)\quad \mbox{so that}
\quad
H(t)=\zg_0\ZSUno (-1)^n\frac{1}{\zl_n} z_n(t) 
\,.
 \]
 
 Now we impose $ \xi_0=0 $ and we choose $ f $ as in~(\ref{eq:expreF}). Then we have
 \begin{multline*}
 \zthe(L,t)=\left (\zg_0^2\ZSUno \frac{(-1)^n}{\zl_n}\right )\intt g(s)\ZD s-\zg_0 \intt g(t-s)H(s)\ZD s\\
 - \zg_0 \intt \beta(s)\left [\int_0^{t-s}g(t-s-r)H(r)\ZD r\right ]\ZD s\,.
  \end{multline*}
  So, $\beta(t)$ is given by the following deconvolution problem (convergence of the numerical series is discussed below)
  
 \begin{multline}\ZLA{eq:PerlaSERIEnumerica}
\intt \left [\int_0^{t-s}g(t-s-r)H(r) \ZD r\right ]\beta(s)\ZD s\\
 =\left (\zg_0\ZSUno \frac{(-1)^n }{\zl_n}\right )\intt g(s)\ZD s-\zg_0\intt g(t-s)H(s)\ZD s- \zthe(L,t)\\
 =\mbox{known.}  
  \end{multline}
The right hand side is known since $H(s)$ is known from the previous measurement.
 
 Leibniz test shows that the numerical series in~(\ref{eq:PerlaSERIEnumerica}) converges but it converges slowly. In spite of this, the identification procedure can be performed since the sum of the series can be explicitly computed,
 \[
 \ZSUno \frac{(-1)^n }{\zl_n}=\ZSUno  \frac{(-1)^n }{n+(1/2)}=   \frac{\pi}{2}-2\,.
 \]
 This is easily seen by noting that for $x\in (0,1)$ we have
 \begin{align*}
 2\left (\sqrt x-\arctan\sqrt x\right )=\int_0^x \sqrt s\frac{1}{1+s }\ZD s=\int_0^x\left (\sum_{k=0}^{+\ZIN}
 (-1)^k x^{k+(1/2)}
 \right )\ZD x\\
 =\sum_{k=0}^{+\ZIN}(-1)^k\frac{1}{k+1+(1/2)} x^{k+1+(1/2)}=
-\sqrt x \sum_{n=1}^{+\ZIN}(-1)^n\frac{1}{n+(1/2)} x^{n }\,.
 \end{align*} 
  The series converges also for $x=1$ and so Abel  Theorem shows that it converges to a function which is continuous also for $x=1$. The sum of the series is obtained by computing both the sides for $x=1$.

\enddocument